\documentclass{amsart}
\usepackage{amsfonts}
\usepackage{amsmath}
\usepackage{amssymb}
\usepackage[utf8]{inputenc}
\usepackage{mathrsfs}
\newtheorem{theorem}{Theorem}
\newtheorem{corollary}{Corollary}
\newtheorem{lemma}{Lemma}
\newtheorem{proposition}{Proposition}

\newtheorem{definition}{Definition}
\newtheorem{remark}{Remark}

\begin{document}

\title[Graded Isomorphisms on Upper Block Triangular Matrix Algebras]{Graded Isomorphisms on Upper Block Triangular Matrix Algebras}

\author{Alex Ramos}
\address{Unidade Acadêmica de Matemática, Universidade Federal de Campina Grande, Campina Grande, PB, 58429-970, Brazil}
\email{almat100@hotmail.com}

\author{Claudemir Fidelis}
\address{Unidade Acadêmica de Matemática, Universidade Federal de Campina Grande, Campina Grande, PB, 58429-970, Brazil}
\email{claudemir@mat.ufcg.edu.br}

\author{Diogo Diniz}
\address{Unidade Acadêmica de Matemática, Universidade Federal de Campina Grande, Campina Grande, PB, 58429-970, Brazil}
\email{diogo@mat.ufcg.edu.br}



\begin{abstract}
We describe the graded isomorphisms of rings of endomorphisms of graded flags over graded division algebras. As a consequence describe the isomorphism classes of upper block triangular matrix algebras (over an algebraically closed field of characteristic zero) graded by a finite abelian group.
\end{abstract}
\maketitle

\section{Introduction}

In this paper we describe the isomorphisms of rings of endomorphisms of a graded flag over a graded division algebra. These rings arise in the classification of gradings on an algebra of upper block triangular matrices, in \cite{VZ} gradings by a finite abelian group were classified under the hypothesis that the base field is algebraically closed of characteristic zero. As it turns out any such graded algebra is isomorphic to the ring of endomorphisms of a graded flag. A grading on an upper block triangular matrix algebra is called elementary if the elementary matrices are homogeneous, any grading of this type is isomorphic to the ring of endomorphisms of a graded flag over a field. Isomorphisms of elementary gradings were described in \cite{BD}, two elementary gradings are isomorphic if and only if the corresponding graded flags are isomorphic up to a shift.

Graded simple algebras (that satisfy the descending chain condition on graded left ideals) are described as the ring of endomorphisms of a graded vector space over a graded division ring (see \cite{EK}, \cite{NO}). More precisely, given a graded simple algebra $R$ satisfying the d.c.c. there exists a pair $(D, V)$, where $D$ is graded division ring  and $V$ is a (finite dimensional) graded vector space over $D$, such that $R$ is isomorphic to $\mathrm{End}_D V$.  Isomorphisms of rings of endomorphisms of graded vector spaces are described in terms of isomorphisms of pairs. In our main result (Theorem \ref{main}) we prove the analogous result for rings of endomorphisms of graded flags.

Gradings on upper block triangular matrix are described in \cite{VZ} as a tensor product of a division grading on a matrix algebra and an elementary grading on an upper block triangular matrix algebra. As a consequence of our main result we determine (Corollary \ref{appl}), in terms of this decompositions, when two gradings are isomorphic. Two gradings are isomorphic if and only if the graded division components and the elementary graded components are isomorphic. The description in \cite{VZ} is for gradings by finite abelian groups and algebras over an algebraically closed field of characteristic zero, it was conjectured in \cite{VZ} that this results hold in general.

\section{Preliminaries}

Let $\mathbb{K}$ a field, $V$ a vector space (over $\mathbb{K}$), a grading by the group $G$ on $V$ (or a $G$-grading on $V$) is a vector space decomposition $V=\oplus_{g\in G}V_g$. The support of the grading is the set $\mathrm{supp}\ V=\{g\in G \mid V_g\neq 0\}$. An element $v\in V$ is homogeneous if $v\in V_g$ for some $g\in G$, if $v\neq 0$ we say that $v$ is homogeneous of degree $g$ and denote $\mathrm{deg}_G\ v=g$.  A subspace $W$ of $V$ is homogeneous if $W=\oplus_{g\in G} W\cap V_g$.

A morphism from $V$ to the $G$-graded vector space $W$ is a linear transformation $f:V\rightarrow W$ such that $f(V_g)\subseteq W_{g}$ for every $g\in G$. If $f$ is an isomorphism of vector spaces then $f^{-1}:W\rightarrow V$ is also a morphism, in this case we say that the graded vector spaces $V$ and $W$ are isomorphic and that $f$ is an isomorphism of graded vector spaces. Let $V$ and $W$ be vector spaces graded by the groups $G$ and $H$ respectively. A linear map $f:V\rightarrow W$ is graded if for every $g\in G$ there exists an $h\in H$ such that $f(V_g)\subseteq W_h$. We say that $f$ is an equivalence of $V$ and $W$ if $f^{-1}$ is also a graded linear transformation. In this case we obtain a bijection $\alpha:\mathrm{supp}\ V \rightarrow \mathrm{supp}\ V^{\prime}$ where $f(V_g)\subseteq W_{\alpha(g)}$.

If $A$ is a $\mathbb{K}$-algebra a $G$-grading on $A$ is a $G$-grading  $A=\oplus_{g\in G} A_g$ on the underlying vector space such that $A_gA_h\subseteq A_{gh}$ for every $g,h\in G$. A subalgebra $B$ of $A$ is homogeneous if $B$ is homogeneous as a subspace of $A$, analogously an ideal of $A$ is homogeneous if the underlying subspace is a homogeneous subspace. Let $A$ and $B$ be algebras graded by the group $G$, a homomorphism of graded algebras is a homomorphism of algebras $\varphi:A\rightarrow B$ such that the underlying linear map is a morphism of graded vector spaces. A grading on the left (resp. right) $A$-module $V$ by the group $G$ is a grading  $V=\oplus_{g\in G} V_g$ on the vector space $V$ such that $A_g V_h\subseteq V_{gh}$ (resp. $V_hA_g\subseteq V_{hg}$) for every $g,h \in G$. We use the convention that morphisms of left modules are written on the right and morphisms of right modules are written on the left. Let $V, V^{\prime}$ be graded $A$-modules, a homomorphism $\psi:V\rightarrow V^{\prime}$ is homogeneous of degree $\tau$ if $\psi(V_g)\subseteq V_{\tau g}^{\prime}$ for every $g\in G$.  

An algebra $D$ graded by $G$ is a graded division algebra if every non-zero homogeneous element is invertible. A graded right module over $D$ is free, the standard results for vector spaces over division algebras hold (see \cite[Proposition 2.5]{TA}). We will refer to a graded right module over a graded division algebra as vector spaces. Let $D$ be a graded division algebra and $V$ a graded vector space over $D$. Denote $\mathrm{End}_D V$ the ring of endomorphisms of the $D$-module $V$. For any $\tau \in G$ the set $(\mathrm{End}_D V)_{\tau}$ of homogeneous endomorphisms of degree $\tau$ is a subspace of $\mathrm{End}_D V$, the sum of such subspaces is direct, moreover if $V$ is of finite dimension we have the equality $\mathrm{End}_D V=\oplus_{\tau \in G}(\mathrm{End}_D V)_{\tau}$ (see \cite{NO}[Corollary I.2.11, pp. 10]), moreover $(\mathrm{End}_D V)_{\tau}(\mathrm{End}_D V)_{\lambda}\subseteq (\mathrm{End}_D V)_{\tau \lambda}$, therefore this decomposition is a grading by the group $G$ on $\mathrm{End}_D V$.

\begin{definition}
Let $D$ be a graded division algebra and $V$ a graded vector space over $D$. A graded flag on $V$ of length $r$ is a chain of homogeneous $D$-subspaces $\mathscr{F}:V_0\subset V_1\subset \cdots \subset V_r$, where $V_0=0$ and $V_r=V$. 
\end{definition}

The set $\mathrm{End}_D\ \mathscr{F}$ of endomorphisms of the flag $\mathscr{F}$ is a subalgebra of $\mathrm{End}_D\ V$. In the next proposition we prove that this is a homogeneous subalgebra.

\begin{proposition}
Let $V$ be a finite dimensional vector space over a graded division algebra $D$. If $\mathscr{F}$ is a graded flag on $V$ then the algebra $\mathrm{End}_D\ \mathscr{F}$ of endomorphisms of the flag is a homogeneous subalgebra of $\mathrm{End}_D\ V$.
\end{proposition}
\textit{Proof.}
Let $\psi$ be an endomorphism of $\mathscr{F}$. We write $\psi=\psi_1+\cdots + \psi_k$, where $\psi_1,\dots, \psi_k$ are homogeneous elements of $\mathrm{End}_D\ V$ of pairwise distinct degrees $\tau_1,\dots, \tau_k$ respectively. Given $v\in V_i\cap V_g$ we have $\psi(v)=\psi_1(v)+\cdots \psi_k(v)$. Since $\psi(v)\in V_i^{\prime}$ and $V_i^{\prime}$ is a homogeneous subspace we conclude that $\psi_j(v)\in V_i^{\prime}$ for $j=1,\dots, k$. This proves that $\psi_j(V_i\cap V_g)\subset V_i^{\prime}$, since $V_i=\oplus_{g\in G}V_i\cap V_g$ we conclude that $\psi_j(V_i)\subseteq V_i^{\prime}$ for $i=1,\dots, r$. This means that $\psi_j$ lies in $\mathrm{End}_D\ \mathscr{F}$ for $j=1,\dots, k$.
\hfill $\Box$


\begin{definition}
Let $D$, $D^{\prime}$ be graded division algebras and $\mathscr{F}:V_0\subset V_1\subset \cdots \subset V_r$, $\mathscr{F}^{\prime}:V_0^{\prime}\subset V_1^{\prime}\subset \cdots \subset V_r^{\prime}$ graded flags of the same length on the vector spaces $V$ and $V^{\prime}$ over $D$ and $D^{\prime}$ respectively. An isomorphism from the pair $(D,\mathscr{F})$ to $(D^{\prime},\mathscr{F}^{\prime})$ is a pair $(\psi_0,\psi_1)$ where $\psi_0:D\rightarrow D^{\prime}$ is an isomorphism of graded algebras and $\psi_1:V\rightarrow V^{\prime}$ is an isomorphism of graded vector spaces such that $\psi_1(V_i)= V_i^{\prime}$, for $i=0,\dots,r$ and  $\psi_{1}(vd)=\psi_1(v)\psi_0(d)$ for every $v \in V$ and every $d\in D$.
\end{definition}

In the next proposition we prove that an isomorphism of pairs $(D,\mathscr{F}) \rightarrow (D^{\prime},\mathscr{F}^{\prime})$ induces an isomorphism from  $\mathrm{End}_D \mathscr{F}$ to $\mathrm{End}_{D^{\prime}} \mathscr{F}^{\prime}$. 

\begin{lemma}\label{endom}
If $R=\mathrm{End}_D \mathscr{F}$ is the ring of endomorphisms of a flag $\mathscr{F}$ over a graded division algebra $D$ then $\mathrm{End}_R\ V_i=D$.
\end{lemma}

\textit{Proof.}
Denote $V_{0}=0$ and let $v\in V_{i}\setminus V_{i-1}$ and $f\in \mathrm{End}_R\ V_i$. If $v$ and $vf$ are linearly independent over $D$ then there exists $r\in R$ such that $rv=0$ and $r(vf)\neq 0$, this is a contradiction because $r(vf)=(rv)f=0$. Hence there exists a $d\in D$ such that $vf=vd$. Given $w\in V_i$ there exists $r\in R$ such that $rv=w$, therefore $wf=(rv)f=r(vd)=wd$.
\hfill $\Box$

\begin{proposition}\label{isopair}
Given an isomorphism $(\psi_0, \psi_1)$ from $(D, \mathscr{F})$ to $(D^{\prime}, \mathscr{F}^{\prime})$, there exists a unique isomorphism of graded algebras $\psi:R\rightarrow R^{\prime}$, where $R=\mathrm{End}_D \mathscr{F}$, $R^{\prime}=\mathrm{End}_{D^{\prime}} \mathscr{F}^{\prime}$, such that $\psi_1(rv)=\psi(r)\psi_1(v)$ for every $r\in R$ and every $v\in V$. Two isomorphisms $(\psi_0, \psi_1)$ and $(\psi_0^{\prime}, \psi_1^{\prime})$ determine the same isomorphism $R\rightarrow R^{\prime}$ of graded algebras if and only if there exists a non-zero homogeneous $d$ in $D^{\prime}_{\epsilon}$, where $\epsilon$ is the neutral element in $G$, such that $\psi_0^{\prime}(x)=d^{-1}\psi_0(x)d$ and $\psi_1^{\prime}(v)=\psi_1(v)d$.
\end{proposition}
\textit{Proof.}
Given $r\in R$ then the mapping $\psi(r)$ given by $v^{\prime}\mapsto \psi_1(r(\psi_1^{-1}(v^{\prime})))$ lies in $R^{\prime}$, moreover $\psi:R\rightarrow R^{\prime}$ is an isomorphism of rings such that $\psi_1(rv)=\psi(r)\psi_1(v)$. It is clear that $\psi$ is unique satisfying this equality and that $\psi$ is an isomorphism of graded rings.

If $(\psi_0, \psi_1)$ and $(\psi_0^{\prime}, \psi_1^{\prime})$ determine the same isomorphism then $(\psi_1)^{-1}\circ \psi_1^{\prime}$ lies in $(\mathrm{End}_R\ V)_{\epsilon}$, therefore Lemma \ref{endom} implies that there exists $d_0\in D_{\epsilon}$ such that $(\psi_1)^{-1}\circ \psi_1^{\prime}(v)=vd_0$ for every $v\in V$. Hence we conclude that $\psi_1^{\prime}(v)=\psi_1(v)d$, for every $v\in V$, where $d=\psi_0(d_0)$. Let $v\in V \setminus \{0\}$ and $x\in D$, then $\psi_1^{\prime}(vx)=\psi_1^{\prime}(v)\psi_0^{\prime}(x)$, therefore $\psi_1(v)\psi_0(x)d=\psi_1(v)d\psi_0^{\prime}(x)$. This implies that $\psi_0^{\prime}(x)=d^{-1}\psi_0(x)d$.
The converse is clear.
\hfill $\Box$

The main result of the paper is Theorem \ref{main} which is the converse of the previous proposition. We also prove an analogous result (Corollary \ref{equiv}) for equivalence of rings of endomorphisms of graded flags.  Let $V$, $W$ be vector spaces graded by the groups $G$ and $H$, respectively, an equivalence from $V$ to $W$ is a linear isomorphism $f:V\rightarrow W$ such that for every $g\in G$ there exists an $h\in H$ such that $f(V_g)=W_h$. Analogously let $A$ and $B$ be algebras graded by the groups $G$ and $H$, an equivalence from $A$ to $B$ is an isomorphism of algebras $\psi:A\rightarrow B$ such that for every $g\in G$ there exists an $h\in H$ such that $\psi(A_g)=B_h$.

\begin{definition}
Let $D$, $D^{\prime}$ be graded division algebras graded by the groups $G$ and $H$ respectively and $V$, $V^{\prime}$ be graded vector spaces over $D$ and $D^{\prime}$ respecively. An equivalence from $(D, V)$ to $(D^{\prime}, V^{\prime})$ is a pair $(\psi_0, \psi_1)$ where $\psi_0:D\rightarrow D^{\prime}$ is an equivalence of graded algebras and $\psi_1:V\rightarrow V^{\prime}$ is an equivalence of vector spaces such that $\psi_1(vd)=\psi_1(v)\psi_0(d)$ for every $v\in V$ and every $d\in D$.
\end{definition}

The homomorphism of algebras determined by an equivalence of pairs may not be an equivalence of graded algebras (see \cite{EK}[pp. 42]). 

We finalize this section with the notion of shift of a grading. Let $V=\oplus_{g\in G} V_g$ be a grading on the vector space $V$ and $g$ an element of $G$. The right shift by $g$ is the grading, denoted $V^{[g]}$, such that $(V^{[g]})_{hg}:=V_{h}$, analogously one defines the left shift by $g$ which is denoted $^{[g]}V$. Let $D$ be a graded algebra and $V$ a graded right $D$-module. Given $g\in G$ then $^{[g^{-1}]} D^{[g]}$ is also a graded algebra and $V^{[g]}$ is a graded right $^{[g^{-1}]}D^{[g]}$-module. Let $\mathscr{F}:V_0\subset V_1\subset \cdots \subset V_r$ be a graded flag on $V$, then $V_0^{[g]}\subset V_1^{[g]}\subset \cdots \subset V_r^{[g]}$ is a graded flag on $V^{[g]}$ that we denote $\mathscr{F}^{[g]}$.

\section{Graded Isomorphisms of Rings of Endomorphisms of Graded Flags}

Let $R=\mathrm{End}_D \mathscr{F}$ and $R^{\prime}=\mathrm{End}_{D^{\prime}} \mathscr{F^{\prime}}$ where $\mathscr{F}$ and $\mathscr{F}^{\prime}$ are flags over the graded division algebras $D$ and $D^{\prime}$ respectively. If there exists a $g\in G$ such that $(^{[g^{-1}]}D^{[g]}, \mathscr{F}^{[g]})$ is isomorphic to $(D^{\prime}, \mathscr{F}^{\prime})$ then Proposition \ref{isopair} implies that $R$ and $R^{\prime}$ are isomorphic as graded rings. In this section, Theorem \ref{main}, we prove the converse of this statement.

\begin{lemma}\label{subm}
	Let $R=\mathrm{End}_D \mathscr{F}$, where  $D$ be a graded division algebra, $V$ a graded right $D$-module and $\mathscr{F}:V_0\subset V_1\subset\dots\subset V_r=V$ is a graded flag on $V$. Then the submodules of $V$ as a left $R$-module are $V_0,V_1,\dots,V_r$.
\end{lemma}
\textit{Proof.}
First, note that  $V_0,V_1,\dots,V_r$ are graded $R$-submodules of $V$. Now, suppose that $W$ is a non-zero graded $R$-submodule of $V$, let $p$ be the integer in $\{1,\dots,r\}$ such that $W\subseteq V_{p}$ and $W\nsubseteq V_{p-1}$. Let $w \in W\setminus V_{p-1}$, there exists a basis  $\beta=\{v_{1},\dots,v_{n}\}$ of $\mathscr{F}$ such that $v_{n_{p-1}+1}=w$, where $n=\mathrm{dim}_D V$, $n_{p-1}=\mathrm{dim}_D V_{p-1}$. Given $v\in V_p$ let $r$ be the $D$-linear endomorphism of $V$ such that $rw=v$ and $rv_i=0$ if $i\neq n_{p-1}+1$, it is clear that $r\in R$. Since $W$ is an $R$-submodule we conclude that $v=rw \in W$, therefore $V_p\subseteq W$.
\hfill $\Box$

\begin{lemma}\label{decomp}
	Let $R=\mathrm{End}_{D} \mathscr{F}$, where $\mathscr{F}$ is a flag over the graded division algebra $D$. If $e$ is an idempotent of $R$ such that $e(V)=V_1$ then $R_1=Re$ is a subring of $R$, the mapping $r\mapsto r\mid_{V_1}$ is an isomorphism from $R_1$ onto $\mathrm{End}_{D} V_1$ and $R=R_1\oplus I_1$, where $I_1=R(1-e)$.
\end{lemma}
\textit{Proof.}
Clearly $R_1$, $I_1$ are left ideals, in particular $R_1$ is a subring. Denote $i:V_1\rightarrow V$ the inclusion, the mapping $r_1\mapsto ir_1e$ is an injective homomorphism from $\mathrm{End}_{D} V_1$ to $R_1$. Moreover given $r\in R_1$ then the restriction $r\mid_{V_1}$ of $r$ to $V_1$ lies in $\mathrm{End}_{D} V_1$ and $i(r\mid_{V_1})e=re=r$, therefore $r_1\mapsto ir_1e$ is an isomorphism of graded rings.
\hfill $\Box$

\begin{remark}\label{rem}
Given $v\in V$ and $r\in R$ the element $re(v)$ lies in $V_1$, therefore $ere(v)=re(v)$. This implies that $ere=re$ for every $r\in R$, therefore $(1-e)r=(1-e)r(1-e)$, hence $I_1$ is also a right ideal. In fact we have $I_1=\mathrm{Ann}_R V_1$, indeed $V_1=e(V)$ implies that $rV_1=0$ if and only if $re=0$. Since   $r=re+r(1-e)$ and $(1-e)e=0$ we conclude that  $re=0$ if and only if $r\in I_1$.
\end{remark}

\begin{lemma}\label{e}
Let $\psi:R\rightarrow R^{\prime}$ be an isomorphism of graded algebras, where $R=\mathrm{End}_{D} \mathscr{F}$ and $R^{\prime}=\mathrm{End}_{D^{\prime}} \mathscr{F}^{\prime}$ are the rings of endomorphisms of the graded flags $\mathscr{F}$ and $\mathscr{F}^{\prime}$ over the graded division algebras $D$ and $D^{\prime}$ respectively.	If $e$ is a non-zero idempotent of $R$ such that $e(V)=V_1$ then $ V_1^{\prime}= \psi(e)(V^{\prime})$.
\end{lemma}
\textit{Proof.}
  The equality $ere=re$ in Remark \ref{rem} implies that $\psi(e)(V^{\prime})$ is a submodule of $V^{\prime}$. It is clear that this is a non-zero submodule, therefore Lemma \ref{subm} implies that $V_1^{\prime}\subset\psi(e)(V^{\prime})$. Let $f^{\prime}$ be a projection of $V^{\prime}$ onto $V_1^{\prime}$ and $f=\psi^{-1}(f^{\prime})$. It follows from the inclusion $ V_1^{\prime}\subseteq \psi(e)(V^{\prime})$ that $\psi(e)\psi(f)=\psi(f)$, therefore $ef=f$. In this case $0\neq f(V)=ef(V)\subseteq V_1$, Lemma \ref{subm} implies that $f(V)=V_1=e(V)$. This implies that $fe=e$, therefore $\psi(e)(V^{\prime})\subseteq \psi(f)(V^{\prime})=V_1^{\prime}$.
\hfill $\Box$

\begin{lemma}\label{basis}
Let $R=\mathrm{End}_{D} \mathscr{F}$ and $R^{\prime}=\mathrm{End}_{D^{\prime}} \mathscr{F}^{\prime}$ be the rings of endomorphisms of the graded flags $\mathscr{F}$ and $\mathscr{F}^{\prime}$ over the graded division algebras $D$ and $D^{\prime}$ respectively. Let $\beta=\{v_1,\dots, v_n\}$ be a basis of $\mathscr{F}$ and let $E_{ij}$ be the $D$-linear endomorphism of $V$ such that $E_{ij}v_k=\delta_{jk}v_i$. If $\psi:R\rightarrow R^{\prime}$ is an isomorphism of graded algebras then $\mathrm{dim}_D V=\mathrm{dim}_{D^{\prime}} V^{\prime}$. Moreover given a non-zero $v_1^{\prime}$ in $\psi(E_{11})(V)$ there exists $v_2^{\prime},\dots, v_n^{\prime}$ in $V^{\prime}$ such that $\beta^{\prime}=\{v_1^{\prime},v_2^{\prime},\dots, v_n^{\prime}\}$ is a basis for $V^{\prime}$ and $\psi(E_{ij})(v_k^{\prime})=\delta_{jk}v_i^{\prime}$ for any $i, j$ such that $E_{ij}\in R$.
\end{lemma}
\textit{Proof.}
Note that $E_{11},\dots, E_{nn}$ are orthogonal idempotents in $R$ such that $E_{11}+\dots+E_{nn}$ is the idendity of $R$. This implies that $V^{\prime}=\oplus_{j=1}^nQ_j^{\prime}$, where $Q_j^{\prime} = \psi(E_{jj})(V^{\prime})$. Since $Q_j^{\prime}$ is a non-zero $D^{\prime}$-subspace of $V^{\prime}$ we conclude that $$\mathrm{dim}_{D^{\prime}}V^{\prime}=\sum_{j=1}^n \mathrm{dim}_{D^{\prime}}Q_j^{\prime}\geq n=\mathrm{dim}_D V.$$ An analogous argument proves the inequality $\mathrm{dim}_D V\geq \mathrm{dim}_{D^{\prime}}V^{\prime}$, therefore we conclude that $\mathrm{dim}_{D^{\prime}}V^{\prime}=\mathrm{dim}_D V$. This implies that $\mathrm{dim}_{D^{\prime}}Q_j^{\prime}=1$ for $j=1,\dots,n$. Given $E_{ij}$ in $R$ we have $\psi(E_{ij})Q_k^{\prime}=\psi(E_{ij})\psi(E_{kk})(V^{\prime})=0$ for any $k\neq j$. If $\psi(E_{ij})Q_{j}^{\prime}= 0$ then $\psi(E_{ij})=0$, which is a contradiction, therefore $\psi(E_{ij})Q_{j}^{\prime}\neq 0$. The equality $E_{ii}E_{ij}=E_{ij}$ implies that  $\psi(E_{ij})Q_{j}^{\prime}\subseteq Q_{i}^{\prime}$, since $Q_{i}^{\prime}$ and $Q_{j}^{\prime}$ are $D^{\prime}$-modules of dimension $1$ we conclude that $v^{\prime}\mapsto \psi(E_{ij})v^{\prime}$ is an isomorphism of $D^{\prime}$-modules from $Q_{j}^{\prime}$ onto $Q_{i}^{\prime}$. Let $v_i^{\prime}$ be the unique element of $Q_i^{\prime}$ such that $\psi(E_{1i}) (v_i^{\prime})=v_1^{\prime}$, $i=2,\dots, n$. Then $\beta^{\prime}=\{v_1^{\prime},\dots, v_n^{\prime}\}$ is a basis of $V^{\prime}$. For any $E_{ij}\in R$ and $k\neq j$ we have $v_k^{\prime}=\psi(E_{ij})(v^{\prime})$ for some $v^{\prime}\in V^{\prime}$, therefore $\psi(E_{ij})(v_k^{\prime})=\psi(E_{ij})\psi(E_{kk})(v^{\prime})=0$. Moreover $E_{1i}\in R$ and $$\psi(E_{1i})\psi(E_{ij})(v_j^{\prime})=\psi(E_{1j})(v_j^{\prime})=v_1^{\prime},$$ therefore $\psi(E_{ij})(v_j^{\prime})=v_i^{\prime}$.
\hfill $\Box$

\begin{remark}\label{gisom}
Let $R$ be a graded simple algebra (i.e. $R^{2}\neq 0$ and $0$, $R$ are the only homogeneous  ideals) that has a minimal graded left ideal $I$. If $V$ is any graded simple left $R$-module, then
there exists $g\in G$ such that $V$ is isomorphic to $I^{[g]}$ as a graded $R$-module (see \cite[Lemma 2.7]{EK}).
\end{remark}

\begin{theorem}\label{main}
Let $D$, $D^{\prime}$ be graded division algebras, $V$, $V^{\prime}$ be graded right modules over $D$ and $D^{\prime}$ respectively and $R=\mathrm{End}_D \mathscr{F}$, $R^{\prime}= \mathrm{End}_{D^{\prime}} \mathscr{F}^{\prime}$, where $\mathscr{F}:V_1\subset \dots \subset V_r$, $\mathscr{F}^{\prime}:V_1^{\prime}\subset \dots \subset V_{r^{\prime}}^{\prime}$ are graded flags on $V$ and $V^{\prime}$ respectively. If $\psi:R\rightarrow R^{\prime}$ is an isomorphism  of graded algebras then $r=r^{\prime}$ and there exists a $g\in G$ and an isomorphism $(\psi_0, \psi_1)$ from $(^{[g^{-1}]}D^{[g]}, \mathscr{F}^{[g]})$ to $(D^{\prime}, \mathscr{F}^{\prime})$ such that $\psi_1(rv)=\psi(r)\psi_1(v)$ for every $r\in R$ and every $v\in V$.
\end{theorem}
\textit{Proof.}
Let $W_1$ be a subspace of $V$ such that $V=V_1\oplus W_1$ and denote $e$ the corresponding projection onto $V_1$. Then $e$ is an idempotent of $R$ such that $e(V)=V_1$. Lemma \ref{decomp} implies that $R=R_1\oplus I_1$, where $R_1=Re$ and $I_1=R(1-e)$. We consider $V^{\prime}$ as a left $R$-module with the action $rv:=\psi(r)(v)$, Lemma \ref{e} implies that $eV^{\prime}=V_1^{\prime}$. Remark \ref{rem} and Lemma \ref{e} imply that $I_1V_1=0$ and $I_1V_1^{\prime}=0$, moreover $V_1$, $V_1^{\prime}$ are simple $R$-modules, hence we conclude that $V_1$ and $V_1^{\prime}$ are simple $R_1$-modules. Remark \ref{gisom} implies that there exists a $g\in G$ and an isomorphism $\psi_1^{\prime}:V_1^{[g]}\rightarrow V_1^{\prime}$ of $R_1$-modules, the equalities $I_1V_1^{\prime}=0$ and $I_1V_1=0$ imply that $\psi_1^{\prime}$ is an isomorphism of $R$-modules. Define $\psi_0:^{[g^{-1}]}D^{[g]}\rightarrow D^{\prime}$ as $$v^{\prime}\psi_0(d)=\psi_1^{\prime}\left((\psi_1^{\prime})^{-1}(v^{\prime})d\right),$$ for $v^{\prime}\in V^{\prime}$ and $d\in D$. Therefore $(\psi_0, \psi_1^{\prime})$ is an isomorphism from $(^{[g^-1]}D^{[g]}, V_1^{[g]})$ to $(D^{\prime}, V_1^{\prime})$.

Let $\beta=\{v_1,\dots, v_n\}$ be a basis for $\mathscr{F}$ and $\beta^{\prime}=\{v_1^{\prime},\dots, v_n^{\prime}\}$ the basis of $V^{\prime}$ given in Lemma \ref{basis} with $v_1^{\prime}=\psi_1^{\prime}(v_1)$. Denote $\psi_1:V\rightarrow V^{\prime}$ the mapping
\begin{equation}\label{psi1}
\psi_1(v_1d_1+\cdots+ v_nd_n)=v_1^{\prime}\psi_0(d_1)+\cdots+v_n^{\prime}\psi_0(d_n).
\end{equation}
Note that $(\psi_0,\psi_1)$ is an isomorphism from $(^{[g^{-1}]}D^{[g]}, V^{[g]})$ to $(D^{\prime}, V^{\prime})$.  Since $E_{ij}v_k=\delta_{jk}v_i$ and $E_{ij}v_k^{\prime}=\delta_{jk}v_i^{\prime}$ we conclude that given $E_{ij}$ in $R$ and $v_k$ in $\beta$ we have $$\psi_1(E_{ij}v_k)=\delta_{jk}v_i^{\prime}=E_{ij}v_k^{\prime}=E_{ij}\psi_1(v_k).$$ It is clear that $\psi_1(v+w)=\psi_1(v)+\psi_1(w)$ and $\psi_1(vd)=\psi_1(v)\psi_0(d)$ for every $v,w \in V$, $d\in D$, therefore we conclude that
\begin{equation}\label{eij}
\psi_1(E_{ij}v)=E_{ij}\psi_1(v),
\end{equation}
for every $E_{ij}\in R$ and every $v\in V$.

Given $d\in D$ denote $r_d$ the $D$-linear endomorphism of $V$ such that $r_d(v_k)=v_kd$ for $k=1,\dots, n$, clearly $r_d\in R$ for every $d\in D$. We have $$r_dv_1^{\prime}=r_d\psi_1(v_1)=\psi_1(r_dv_1)=\psi_1(v_1d)=v_1^{\prime}\psi_0(d).$$ We claim that $r_dv_k^{\prime}=v_k^{\prime}\psi_0(d)$ for $k=2,\dots, n$. Note that for any $d\in D$ and every $E_{ij}\in R$ we have $E_{ij}r_d=r_dE_{ij}$. Let $d_1^{\prime},\dots, d_n^{\prime}$ be the elements in $D^{\prime}$ such that
$$r_dv_k^{\prime}=v_1^{\prime}d_1^{\prime}+\dots+v_k^{\prime}d_k^{\prime}.$$
If $i\neq k$ then $v_i^{\prime}d_i^{\prime}=E_{ii}(r_dv_k^{\prime})=r_dE_{ii}v_k^{\prime}=0$, hence $d_i^{\prime}=0$. Therefore $r_dv_k^{\prime}=v_k^{\prime}d_k^{\prime}$ and we conclude that $$v_1^{\prime}d_k^{\prime}=E_{1k}v_k^{\prime}d_k^{\prime}=E_{1k}r_dv_k^{\prime}=r_dv_1^{\prime}=v_1^{\prime}\psi_0(d),$$ this implies that $d_k^{\prime}=\psi_0(d)$ and the claim is proved. We have $$\psi_1(r_dv_k)=\psi_1(v_kd)=v_k^{\prime}\psi_0(d)=r_dv_k^{\prime}=r_d\psi_1(v_k),$$ therefore
\begin{equation}\label{rd}
\psi_1(r_dv)=r_d\psi_1(v),
\end{equation}
for every $d\in D$ and every $v\in V$.
Given $r\in R$, let $d_{ij}$, $1\leq i,j \leq n$ be the elements of $D$ such that $rv_k=v_1d_{1k}+\dots+ v_nd_{nk}$. Then $r=\sum r_{d_{ij}}E_{ij}$, where the sum is over the $i,j$ such that $E_{ij}\in R$. Therefore $R$ is generated as a ring by the $E_{ij}$ that lie in $R$ together with the elements $\{r_{d} \mid d\in D\}$. Therefore (\ref{eij}) and (\ref{rd}) imply that $\psi_1(rv)=r\psi_1(v)$ for every $r\in R$ and every $v\in V$. This implies that $\psi_1(V_1)\subset \dots \subset \psi_1(V_r)$ are $R$-submodules of $V^{\prime}$, analogously $\psi_1^{-1}(V_1^{\prime})\subset \dots \subset \psi_1^{-1}(V_{r^{\prime}}^{\prime})$ are $R$-submodules of $V$. Therefore Lemma \ref{subm} implies that $r=r^{\prime}$ and $\psi_1(V_i)=V_i^{\prime}$ for $i=1,\dots, r$. 
\hfill $\Box$

Let $G$ be a group. Given $D$ a graded division algebra and $\textbf{g}\in G^{n}$, where $n$ is a natural number, denote $V(D,n, \textbf{g})$ the graded right $D$-module $\oplus_{i=1}^n\ ^{[g_i]}D$.

\begin{remark}\label{isom}
The pairs $(D,V(D,n, \textup{\textbf{g}}))$ and $(D^{\prime},V(D^{\prime}, n^{\prime},\textup{\textbf{g}}^{\prime}))$ are isomorphic if and only if, $D\cong D^{\prime}$, $n=n^{\prime}$ and there exists $h_1,\dots, h_n\in \mathrm{supp}\ D$, $\sigma \in S_n$ such that $g_i=g_{\sigma(i)}^{\prime}h_{\sigma(i)}$ for $i=1,\dots, n$.
\end{remark}

Given a triple $(D,\textbf{m}, \textbf{g})$, where $D$ is an algebra with a division grading by the group $G$, $\textbf{m}=(m_1,\dots, m_s)$ and $\textbf{g}=(g_1,\dots, g_n)\in G^n$, $n=m_1+\cdots +m_s$, denote $\mathscr{F}(D,\textbf{m}, \textbf{g})$ the graded flag $V_0\subset V_1\subset \cdots \subset V_s$, where $V_0=0$ and $V_i=\oplus_{j=1}^{n_i}\ ^{[g_j]}D$, $n_i=m_1+\cdots  +m_i$, for $i=1,\dots, s$. The ring $R=\mathrm{End}_D \mathscr{F}$ of endomorphisms of this flag is denoted $\mathscr{A}(D,\textbf{m}, \textbf{g})$.

\begin{lemma}
Let   $\mathscr{F}(D,\textup{\textbf{m}}, \textup{\textbf{g}}):V_0\subset \dots \subset V_s$ and $\mathscr{F}(D^{\prime},\textup{\textbf{m}}^{\prime}, \textup{\textbf{g}}^{\prime}):V_0^{\prime}\subset \dots \subset V_s^{\prime}$ be graded flags of the same length and let $\psi_0:D\rightarrow D^{\prime}$ be an isomorphism of graded algebras. There exists an isomorphism $\psi_1:V_s\rightarrow V_s^{\prime}$ of graded vector spaces such that $(\psi_0, \psi_1)$ is an isomorphism of pairs if and only if $\textup{\textbf{m}}=\textup{\textbf{m}}^{\prime}$ and there exists $h_1,\dots, h_n \in \mathrm{supp}\ D$, where $n=n_1+\cdots +n_s$, and $\sigma\in S_{m_1}\times \cdots \times S_{m_s}$ such that $g_i=g_{\sigma(i)}^{\prime}h_{\sigma(i)}$.
\end{lemma}
\textit{Proof}
Denote $\textbf{m}=(m_1,\dots, m_s)$ and $\textbf{m}^{\prime}=(m_1^{\prime},\dots, m_{s^{\prime}}^{\prime})$. Assume that there exists an isomorphism of graded vector spaces $\psi_1:V_s\rightarrow V^{\prime}_s$ such that $(\psi_0, \psi_1)$ is an isomorphism of pairs. In this case the the pairs $(D,V_{i+1}/V_i)$ and $(D^{\prime}, V_{i+1}^{\prime}/V_i^{\prime})$ are isomorphic, $i=0,\dots, s-1$. The vector space $V_{i+1}/V_i$ has a basis of homogeneous elements of degrees $(g_{n_i+1},\dots, g_{n_i+m_{i+1}})$, where $n_0=0$ and $n_i=n_{i-1}+m_i$ for $i=1,\dots,s$. Analogously $V_{i+1}^{\prime}/V_i^{\prime}$ has a basis of homogeneous elements of degrees $(g_{n_i^{\prime}+1}^{\prime},\dots, g_{n_{i}^{\prime}+m_{i+1}^{\prime}})$. Remark \ref{isom} implies $m_i=m_i^{\prime}$ and that there exists $h_{n_i+1},\dots, h_{n_{i}+m_i} \in \mathrm{supp}\ D$ and a permutation $\sigma_i$ of the elements $\{n_i+1,\dots, n_i+m_i\}$ such that $g_j=g_{\sigma(j)}^{\prime}h_{\sigma(j)}$ for every $j\in \{n_i+1,\dots, n_i+m_i\}$. Therefore we conclude that $\textbf{m}=\textbf{m}^{\prime}$, moreover if $\sigma = \sigma_1\cdots \sigma_s$ then $g_i=g_{\sigma(i)}^{\prime}h_{\sigma(i)}$ for $i=1,\dots,n$. To prove the converse let $\{v_1^{\prime},\dots, v_n^{\prime}\}$ be a basis of $\mathscr{F}(D,\textbf{m}^{\prime}, \textbf{g}^{\prime})$ of homogeneous elements with $\mathrm{deg}_G v_i^{\prime}=g_i^{\prime}$. Let $d_1,\dots, d_n$ be non-zero homogeneous elements of degrees $h_1,\dots, h_n$ respectively. Now let $w_i=v_{\sigma(i)}^{\prime}d_{\sigma(i)}$, the set $\{w_1,\dots, w_n\}$ is a basis of $\mathscr{F}^{\prime}$ and $\mathrm{deg}_G w_i=g_i$ for $i=1,\dots, n$. Now let $\{v_1,\dots, v_n\}$ be a basis for $\mathscr{F}$ of homogeneous elements of degrees $(g_1,\dots, g_n)$. Then the map $\psi_1(v_1d_1+\dots+v_nd_n)= w_1\psi_0(d_1)+\dots + w_n \psi_0(d_n)$ is an isomorphism of graded vector spaces such that $(\psi_0, \psi_1)$ is an isomorphism of pairs.
\hfill $\Box$

\begin{corollary}\label{cisom}
The algebras $\mathscr{A}(D,\textup{\textbf{m}}, \textup{\textbf{g}})$ and $\mathscr{A}(D^{\prime},\textup{\textbf{m}}^{\prime}, \textup{\textbf{g}}^{\prime})$ are isomorphic if and only if $\textup{\textbf{m}}=\textup{\textbf{m}}^{\prime}$ and there exists a $g\in G$ such  $^{[g^{-1}]}D^{[g]}\cong D^{\prime}$ and there exists $h_1,\dots, h_n\in \mathrm{supp}\ D$ and $\sigma\in S_{m_1}\times \cdots \times S_{m_s}$ such that $g_i^{\prime}=g_{\sigma(i)}h_{\sigma(i)}g$ for $i=1,\dots,n$.
\end{corollary}
\textit{Proof.}
Theorem \ref{main} implies that the algebras $\mathscr{A}(D,\textbf{m}, \textbf{g})$ and $\mathscr{A}(D^{\prime},\textbf{m}^{\prime}, \textbf{g}^{\prime})$ are isomorphic  if and only if the flags $\mathscr{F}(D,\textbf{m}, \textbf{g})$ and $\mathscr{F}(D^{\prime},\textbf{m}^{\prime}, \textbf{g}^{\prime})$ are of the same length and there exists a $g\in G$ and an isomorphism $(\psi_0, \psi_1)$ from $(D^{\prime}, \mathscr{F}^{\prime})$ to $(^{[g^{-1}]}D^{[g]}, \mathscr{F}^{[g]})$. Note that $\mathscr{F}(D,\textbf{m}, \textbf{g})^{[g]}=\mathscr{F}(D,\textbf{m}, \textbf{k})$, where $\textbf{g}=(g_1,\dots, g_s)$ and $\textbf{k}=(g_1g,\dots, g_sg)$, moreover $\mathrm{supp}\ ^{[g^{-1}]}D^{[g]}=g^{-1}(\mathrm{supp}\ D) g$. The result now follows from the previous lemma.
\hfill $\Box$

Now we consider equivalence between rings of endomorphisms of graded flags. In the case of rings of endomorphisms of vector spaces we have the following result.

\begin{proposition}\cite{EK}[Proposition 2.33]\label{gsimpleequiv}
Let $D$, $D^{\prime}$ be graded division algebras, $V$, $V^{\prime}$ be graded right modules over $D$ and $D^{\prime}$ respectively. If $\psi:R\rightarrow R^{\prime}$ is an equivalence  of graded algebras, where $R=\mathrm{End}_D V$, $R^{\prime}= \mathrm{End}_{D^{\prime}} V^{\prime}$ then there exists an equivalence of pairs $(\psi_0, \psi_1)$ from $(D, V)$ to $(D^{\prime}, V^{\prime})$ such that $\psi_1(rv)=\psi(r)\psi_1(v)$ for every $r\in R$ and every $v\in V$.
\end{proposition}

As in Theorem \ref{main} we obtain the analogous result for rings of endomorphisms of graded flags. We remark that the converse of this proposition does not hold.

\begin{corollary}\label{equiv}
Let $D$, $D^{\prime}$ be graded division algebras, $V$, $V^{\prime}$ be graded right modules over $D$ and $D^{\prime}$ respectively and $R=\mathrm{End}_D \mathscr{F}$, $R^{\prime}= \mathrm{End}_{D^{\prime}} \mathscr{F}^{\prime}$, where $\mathscr{F}$, $\mathscr{F}^{\prime}$ are graded flags on $V$ and $V^{\prime}$ respectively. If $\psi:R\rightarrow R^{\prime}$ is an equivalence  of graded algebras then there exists an equivalence of pairs $(\psi_0, \psi_1)$ from $(D, \mathscr{F})$ to $(D^{\prime}, \mathscr{F}^{\prime})$ such that $\psi_1(rv)=\psi(r)\psi_1(v)$ for every $r\in R$ and every $v\in V$.
\end{corollary}
\textit{Proof.}
Let $e$ be a projection of $V$ onto $V_1$ and $R_1=Re$. It follows from Lemma \ref{e} that $e^{\prime}=\psi(e)$ is a projection onto $V^{\prime}_1$, let $R_{1}^{\prime}=R^{\prime}e^{\prime}$. Lemma \ref{decomp} implies that $r\mapsto r\mid_{V_1}$ is an isomorphism from $R_1$ onto $\mathrm{End}_D V_1$ and $r^{\prime}\mapsto r^{\prime}\mid_{V_1^{\prime}}$ is an isomorphism from $R_1^{\prime}$ onto $\mathrm{End}_D V_1^{\prime}$. The restriction of $\psi$ to $R_1$ is an equivalence from $R_1$ to $R_{1}^{\prime}$, hence Proposition \ref{gsimpleequiv} implies that there exists an equivalence of pairs $(\psi_0, \psi_1)$ from $(D, V_1)$ to $(D^{\prime}, V^{\prime}_1)$ such that $\psi_1(rv)=\psi(r)\psi_1(v)$ for every $r\in R_1$ and every $v\in V_1$. Lemma \ref{decomp} implies that $\psi_1(rv)=\psi(r)\psi_1(v)$ for every $r\in R$. Now let $\beta=\{v_1,\dots, v_n\}$ be a basis of $\mathscr{F}$ and $\beta^{\prime}=\{v_1^{\prime},\dots, v_n^{\prime}\}$ be the basis obtained from Lemma \ref{basis} with $v_1^{\prime}=\psi_1(v_1)$. Then the map $\psi_1:V\rightarrow V^{\prime}$ given by $\psi_1(v_1d_1+\cdots + v_nd_n)=v_1^{\prime}\psi_0(d_1)+\cdots+v_n^{\prime}\psi_0(d_n)$ is an equivalence of graded vector spaces such that $(\psi_0, \psi_1)$ is an equivalence of pairs from $(D, \mathscr{F})$ to $(D^{\prime}, \mathscr{F}^{\prime})$. It follows from the proof of Theorem \ref{main} that $\psi_1(rv)=\psi(r)\psi_1(v)$ for every $r\in R$ and every $v\in V$.
\hfill $\Box$

As a consequence of Corollary \ref{equiv} we obtain the following.

\begin{corollary}\label{cequiv}
If the algebras $\mathscr{A}(D,\textup{\textbf{m}}, \textup{\textbf{g}})$ and $\mathscr{A}(D^{\prime},\textup{\textbf{m}}^{\prime}, \textup{\textbf{h}})$ are equivalent then $D$ is equivalent to $D^{\prime}$, $\textup{\textbf{m}}=\textup{\textbf{m}}^{\prime}$ and there exists a bijection $\lambda$ from $\{g(\mathrm{supp}\ D)\mid g\in \mathrm{supp}\ V\}$ to $\{g^{\prime}(\mathrm{supp}\ D^{\prime})\mid g^{\prime}\in \mathrm{supp}\ V^{\prime}\}$ and a permutation $\sigma\in S_{m_1}\times \cdots \times S_{m_s}$ such that $h_{\sigma(i)}\mathrm{supp}\ D^{\prime}=\lambda(g_{i}\mathrm{supp}\ D)$ for $i=1,\dots,n$.
\end{corollary}
\textit{Proof.}
Corollary \ref{equiv} implies that there exists an equivalence from $\mathscr{F}(D,\textbf{m}, \textbf{g})$ to $\mathscr{F}(D^{\prime},\textbf{m}^{\prime}, \textbf{g}^{\prime})$. Let $(\psi_0,\psi_1)$ be an equivalence of pairs this implies that $D$ is equivalent to $D^{\prime}$, $\textbf{m}=\textbf{m}^{\prime}$. Given $g\in (\mathrm{supp}\ V)$ there exists a unique $g^{\prime}\in \mathrm{supp}\ V^{\prime}$ such that $\psi_1(V_g)=V_{g^{\prime}}^{\prime}$ hence we have a bijection $g\mapsto g^{\prime}$. Since $\psi_1(vd)=\psi_1(v)\psi_0(d)$ for every $v\in V$ and every $d\in D$ we conclude that $g(\mathrm{supp}\ D)=h(\mathrm{supp}\ D)$ if and only if $g^{\prime}(\mathrm{supp}\ D^{\prime})=h^{\prime}(\mathrm{supp}\ D^{\prime})$, therefore $g(\mathrm{supp}\ D)\mapsto g^{\prime}(\mathrm{supp}\ D^{\prime})$ is a bijection which we denote $\lambda$. Let $\beta=\{v_1,\dots, v_n\}$ be a basis of $\mathscr{F}$ of homogeneous elements of degrees $g_1,\dots, g_n$ respectively. Then $\{v_1^{\prime},\dots, v_n^{\prime}\}$, where $v_i^{\prime}=\psi_1(v_i)$ is a basis of $\mathscr{F}^{\prime}$ of homogeneous elements of degrees $g_1^{\prime},\dots, g_n^{\prime}$ respectively. Therefore there exists a permutation $\sigma\in S_{m_1}\times \cdots \times S_{m_s}$ such that $\lambda(g_{i}\mathrm{supp}\ D)=g_{i}^{\prime}\mathrm{supp}\ D^{\prime}=h_{\sigma(i)}\mathrm{supp}\ D^{\prime}$.
\hfill $\Box$

\section{Gradings on Upper Block Triangular Matrix Algebras}

Let $\textbf{p}=(p_1,\dots, p_s)$ be an $s$-tuple of natural numbers. Define inductively $n_0=0$, $n_i=p_1+\cdots +p_i$, for $i=1,\dots,s$ and denote $P_i=\{n_{i-1}+1,\dots, n_{i}\}$,  $i=1,\dots, s$. Given $i\in \{1,\dots, n_s\}$ there exists a unique $k\in \{1,\dots,s\}$ such that $i\in P_k$. If $i\in P_k$ and $j\in P_l$ then the elementary matrix $e_{ij}$, with $1$ as the only non-zero entry in the $i$-th line and $j$-th column, lies in $R=UT(p_1,\dots, p_s)$ if and only if $k\leq l$. The set of elementary matrices $e_{ij}$ that lie in $R$ form a basis for this algebra, we refer to this as the canonical basis of $R$.

\begin{definition}
A grading $R=\oplus_{g\in G}R_g$ by the group $G$ on the algebra $R=UT(p_1,\dots, p_s)$ of upper-block triangular matrix algebras is an elementary grading if every elementary matrix in the canonical basis of $R$ is homogeneous. 
\end{definition}

If $R=UT(p_1,\dots, p_s)$ has an elementary grading then there exists an $n$-tuple $\textbf{g}=(g_1,\dots, g_n)$ of elements of $G$, where $n=p_1+\cdots+p_s$, such that $\mathrm{deg}_G e_{ij}=g_ig_j^{-1}$. Conversely given an $n$-tuple $\textbf{g}=(g_1,\dots, g_n)$ we denote $R_g$ the subspace of $R$ generated by the matrices $e_{ij}$ in the canonical basis of $R$, where $i,j$ are such that $g_ig_j^{-1}=g$. Then $R=\oplus_{g\in G}R_g$ is an elementary grading on $R$. 

Any elementary grading is isomorphic to the algebra of endomorphisms of a suitable flag over a field. In \cite{BD} the authors proved that two such rings of endomorphisms are isomorphic if and only if the graded flags are isomorphic up to a shift, moreover the authors determined in terms of the associated tuples when two elementary gradings are isomorphic. In this section we consider the analogous results for arbitrary gradings on upper-block triangular matrix algebras.

\begin{definition}
A $G$-grading $R=\oplus_{g\in G}R_g$ on the algebra $R$ is fine if $\mathrm{R_g}\leq 1$ for every $g\in G$.
\end{definition}

\begin{lemma}\cite{EK}[Lemma 2.20]\label{fine}
Let $R$ be a matrix algebra over an algebraically closed field and let $R=\oplus_{g\in G} R_g$ be a grading by the group $G$. Then the following conditions are equivalent:
\begin{enumerate}
\item[1)] $dim R_g\leq 1$ for all $g\in G$;
\item[2)] $dim R_e=1$;
\item[3)] $R$ is a graded division algebra.
\end{enumerate}
\end{lemma}

Let $G$ be a group, not necessarily abelian, $R=M_n(\mathbb{K})$ with an elementary grading induced by $(g_1,\dots, g_n)\in G^{n}$ and $D$ an algebra graded by $G$. The tensor product $R\otimes D$ has a $G$-grading such that $\mathrm{deg}_G(e_{ij}\otimes d)=g_i\mathrm{deg}_G d g_j^{-1}$. If the group $G$ is abelian the tensor product on any two $G$-graded algebras $R$ and $S$ has a canonical grading where $(R\otimes S)_g=\oplus_{hk=g}R_{h}\otimes S_{k}$, this coincides with the previous grading if $R$ is a matrix algebra $M_n(\mathbb{K})$ with an elementary grading. The main result of \cite{VZ} is that if $G$ is a finite abelian group and the base field $\mathbb{K}$ is algebraically closed of characteristic zero then every grading on an upper-block triangular matrix algebra is isomorphic as a graded algebra to a tensor product of an upper-block triangular matrix algebra with an elementary grading and full matrix algebra with a fine grading. 

\begin{theorem}\cite[Theorem 3.2]{VZ}\label{VZ}
Let $G$ be a finite abelian group and $UT(d_1,\dots, d_m)$ an upper block triangular matrix algebra over an algebraically closed field $\mathbb{K}$ of characteristic zero. Then there exists a decomposition $d_1 = tp_1, \dots, d_m=tp_m$, a subgroup $H\subset G$, and an $n$-tuple $(g_1, \dots, g_n)\in G^n$ , where $n = p_1 + \cdots + p_m$ such that $UT(d_1,\dots, d_m)$ is isomorphic to $M_t(\mathbb{K})\otimes UT(p_1,\cdots , p_m)$ as a $G$-graded algebra where $M_t(\mathbb{K})$ is an $H$-graded algebra with a “fine” grading with support $H$ and $UT(p_1,\cdots , p_m)$ has an elementary grading defined by $(g_1, \dots, g_n)$.
\end{theorem}

Since $\mathbb{K}$ is algebraically closed Lemma \ref{fine} implies that a fine grading on $M_t(\mathbb{K})$ is a division grading. In the next proposition we prove that the algebra $M_t(\mathbb{K})\otimes UT(p_1,\cdots , p_m)$ in the previous theorem is isomorphic to the ring of endomorphisms of a flag over $D=M_t(\mathbb{K})$.

\begin{proposition}\label{tensor}
Let $G$ be a group, $D$ a graded division algebra, $\textup{\textbf{p}}=(p_1,\dots, p_s)$ an $s$-tuple of natural numbers and $\textup{\textbf{g}}=(g_1,\dots, g_n)$ an $n$-tuple of elements of $G$, where $n=p_1+\cdots +p_s$. The algebra $UT(p_1,\dots, p_s)\otimes_{\mathbb{K}} D$ with the grading such that $\mathrm{deg}_G\ e_{ij}\otimes d =g_i (\mathrm{deg}_G\ d) g_j^{-1}$ is isomorphic to the algebra $\mathscr{A}(D,\textup{\textbf{p}}, \textup{\textbf{g}})$.
\end{proposition}

\textit{Proof.}
Let $n=p_1+\cdots +p_s$, we consider in $M_n(\mathbb{K})\otimes_{\mathbb{K}} D$ the grading such that $\mathrm{deg}_G\ e_{ij}\otimes d =g_i (\mathrm{deg}_G\ d) g_j^{-1}$. The algebra $UT(p_1,\dots, p_s)\otimes_{\mathbb{K}} D$ is a homogeneous subalgebra of $M_n(\mathbb{K})\otimes_{\mathbb{K}} D$.  Note that $\mathscr{A}(D,\textbf{p}, \textbf{g})$ is a homogeneous subalgebra of $\mathrm{End}_D V$, where $V=V(D,n, \textbf{g})$. Let $\beta=\{v_1,\dots, v_n\}$ be a basis of $\mathscr{F}(D,\textbf{p}, \textbf{g})$ of homogeneous elements such that $\mathrm{deg}_G\ v_i=g_i$. Denote $\varphi:M_n\otimes_{\mathbb{K}} D\rightarrow \mathrm{End}_D V$ the isomorphism given by $$\varphi((\lambda_{ij})\otimes d)v_k=\sum_i v_i(\lambda_{ik}d), k=1,\dots,n.$$ We denote $E_{ij}$ the element of $\mathrm{End}_D V$ such that $E_{ij}v_k=\delta_{jk}v_i$ for $1\leq i,j,k \leq n$. Given $d\in D$ we denote $r_d$ the element of $\mathrm{End}_D V$ such that $r_dv_k=v_kd$ for $k=1,\dots,n$. An endomorphism $r\in \mathrm{End}_D V$ is determined by the equalities $rv_j=\sum_i v_id_{ij}$, where $d_{ij}\in D$, $j=1,\dots, n$. We have $r=\sum_{ij}E_{ij}r_{d_{ij}}$, moreover element $r$ lies in $\mathscr{F}(D,\textbf{p}, \textbf{g})$ if and only if $d_{ij}=0$ whenever $i$ and $j$ are such that $i\in P_k$, $j\in P_l$ and $k>l$. This implies that $\mathscr{F}(D,\textbf{p}, \textbf{g})$ is generated by the elements $\{E_{ij}\mid i\in P_k, j\in P_l, 1\leq k\leq l \leq s\}$ together with the elements $\{r_d\mid d\in D\}$. Since $\varphi(e_{ij}\otimes d)=E_{ij}r_d$ we conclude that $\varphi$ maps $UT(p_1,\dots, p_s)\otimes_{\mathbb{K}} D$ onto $\mathscr{A}(D,\textbf{p}, \textbf{g})$.
\hfill $\Box$

\begin{corollary}\label{appl}
Let $S$ and $S^{\prime}$ denote the algebras $M_t(\mathbb{K})\otimes_{\mathbb{K}} UT(p_1,\cdots , p_m)$ and $M_{t^{\prime}}(\mathbb{K})\otimes_{\mathbb{K}} UT(p_1^{\prime},\cdots , p_m^{\prime})$ respectively, where $M_t(\mathbb{K})$, $M_{t^{\prime}}(\mathbb{K})$ are algebras with a fine grading and $UT(p_1,\cdots , p_m)$, $UT(p_1^{\prime},\cdots , p_m^{\prime})$ have elementary gradings defined by the tuples $\textup{\textbf{g}}$, $\textup{\textbf{g}}^{\prime}$ of elements of $G$, respectively. The algebras $S$ and $S^{\prime}$ are isomorphic if and only if $M_t(\mathbb{K})\cong M_{t^{\prime}}(\mathbb{K})$, $(p_1,\cdots , p_m)=(p_1^{\prime},\cdots , p_m^{\prime})$ and there exists a $g\in G$, $h_1,\dots, h_n\in \mathrm{supp}\ M_t$ and $\sigma\in S_{p_1}\times \cdots \times S_{p_m}$ such that $g_i^{\prime}=g_{\sigma(i)}h_{\sigma(i)}g$ for $i=1,\dots,n$.
\end{corollary}
\textit{Proof.}
Proposition \ref{tensor} implies that $S\cong \mathscr{A}(D,\textbf{p},\textbf{g})$, where $D=M_t(\mathbb{K})$, $\textbf{p}=(p_1,\cdots , p_m)$ and $S^{\prime}\cong \mathscr{A}(D^{\prime},\textbf{p}^{\prime},\textbf{g}^{\prime})$, where $D^{\prime}=M_{t^{\prime}}(\mathbb{K})$, $\textbf{p}^{\prime}=(p_1^{\prime},\cdots , p_m^{\prime})$. Since the group is abelian $^{[g^{-1}]} D ^{[g]}=D$ for any $g\in G$, therefore the result follows from Corollary \ref{cisom}.
\hfill $\Box$

The converse of Corollary \ref{equiv} does not hold in general. We prove that equivalence of gradings distinguishes elementary and non-elementary gradings.

\begin{corollary}\label{appleq}
Let $S= UT(p_1,\cdots , p_m)$ be an upper-block triangular matrix algebra with an elementary grading by the group $G$ defined by the tuple $\textup{\textbf{g}}$ and $S^{\prime}= UT(p_1^{\prime},\cdots , p_m^{\prime})$ graded by the group $H$ respectively. The algebras $S$ and $S^{\prime}$ are equivalent if and only if $S^{\prime}$ has an elementary grading $(p_1,\cdots , p_m)=(p_1^{\prime},\cdots , p_m^{\prime})$ and there exists a bijection $\lambda$ from $\{g_1,\dots, g_n\}$ to $\{h_1,\dots, h_n\}$ such that the equality $\lambda(g_i)\lambda(g_j)^{-1}=\lambda(g_k)\lambda(g_l)^{-1}$ holds whenever $g_ig_j^{-1}=g_kg_l^{-1}$  and a permutation $\sigma\in S_{p_1}\times \cdots \times S_{p_m}$ such that $h_{\sigma(i)}=\lambda(g_{i})$ for $i=1,\dots,n$, where $h_i=(\mathrm{deg}_H e_{1i})^{-1}$.
\end{corollary}
\textit{Proof.}
If $S^{\prime}$ is equivalent to $S$ then Corollary \ref{cequiv} and Theorem \ref{VZ} imply that $S^{\prime}$ has an elementary grading. Note that $(h_1,\dots, h_n)$, where $h_i=(\mathrm{deg}_HE_{1i})^{-1}$ induces the elementary grading in $S^{\prime}$. The existence of $\lambda$ and $\sigma$ also follows from Corollary \ref{cequiv}. 
To prove the converse note that $S\cong \mathscr{A}(\mathbb{K},\textbf{p}, \textbf{g})$ and $S^{\prime}\cong \mathscr{A}(\mathbb{K},\textbf{p}, \textbf{h})$. Let $\{v_1,\dots, v_n\}$ be the canonical basis of $\mathscr{F}(\mathbb{K},\textbf{p}, \textbf{g})$ and $\{v_1^{\prime},\dots, v_n^{\prime}\}$ be a basis of $\mathscr{F}(\mathbb{K},\textbf{p}, \textbf{h})$ where $\mathrm{deg}_Hv_i^{\prime}=\lambda(g_i)$. The linear map $\psi_1$ such that $\psi_1(v_i)=v_i^{\prime}$ is an equivalence of vector spaces. Let $\psi:\mathscr{A}(\mathbb{K},\textbf{p}, \textbf{g})\rightarrow \mathscr{A}(\mathbb{K},\textbf{p}, \textbf{h})$ be the homomorphism of algebras such that $\psi_1(rv)=\psi(r)\psi_1(v)$ for every $r\in \mathscr{A}(\mathbb{K},\textbf{p}, \textbf{g})$ and every $v\in V$. Then $\psi(e_{ij})$ is homogeneous of degree $\lambda(g_{i})\lambda(g_j)^{-1}$. This is an equivalence of algebras because $\lambda(g_i)\lambda(g_j)^{-1}=\lambda(g_k)\lambda(g_l)^{-1}$ whenever $g_ig_j^{-1}=g_kg_l^{-1}$.
\hfill $\Box$

\end{document}